\documentclass[11pt]{article}
\usepackage{amsmath,amsthm,latexsym,amsfonts,graphicx,amssymb,epsfig}
\title{Proper actions of lattices on contractible manifolds}
  \author{Mladen Bestvina and Mark Feighn \thanks{The authors
  gratefully acknowledge the support by the National Science
  Foundation.}}  \date{\today}

\newtheorem{thm}{Theorem}[section]
\newtheorem{lemma}[thm]{Lemma}
\newtheorem{cor}[thm]{Corollary}

\theoremstyle{remark}

\newtheorem{remark}{Remark}
\newtheorem{example}[thm]{Example}
\def\bigsum{\sum}
\def\ha{\hat\alpha}
\def\frak{\mathfrak}
  
\def\L{L}
\def\R{{\mathbb R}}
\def\Z{{\mathbb Z}}
\def\C{{\mathbb C}}
\def\H{{\mathbb H}}
\def\G{{\Gamma}}  
\def\O{{\cal O}} 
\def\g{{\frak g}} 
\def\n{{\frak n}} 

\def\Q{{\mathbb Q}}
\def\a{\alpha}
\def\b{\beta}

\def\obdim{\operatorname{obdim}}

\def\cd{\operatorname{cd}}
\def\rank{\operatorname{rank}}
     
\begin{document}
\maketitle
\begin{abstract}
  Every lattice $\Gamma$ in a connected semi-simple Lie group $G$ acts
  properly discontinuously by isometries on the contractible manifold
  $G/K$ ($K$ a maximal compact subgroup of $G$). We prove that if
  $\Gamma$ acts on
  a contractible manifold $W$ and if either 

1) the action
  is properly discontinuous, or 

2) $W$ is equipped with a complete Riemannian metric, the action is by
  isometries and with unbounded orbits, $G$ is simple with finite
  center and rank $>1$,

  then $\dim W\geq\dim G/K$.
\end{abstract}

\section{Introduction}

Let $G$ be a connected semi-simple Lie group and $K$ a maximal compact subgroup
of $G$. It is a theorem of Malcev \cite{malcev:k} and Iwasawa \cite{iwasawa:k}
that the homogeneous space $G/K$ (even without the assumption of
semi-simplicity) is a contractible manifold (see \cite[Theorem
XV.3.1]{hochschild} and also \cite{mostow:sag} for the case when $G$ has
finitely many components). When $G$ has finite center, $G/K$ is a symmetric
space of noncompact type. Every lattice $\G\subset G$ acts properly
discontinuously on $G/K$. The main theorem of this paper is:

\begin{thm} \label{main}
  If $\G$ acts properly discontinuously (by homeomorphisms) on a contractible
  manifold $W$, then $\dim W\geq \dim G/K$.
\end{thm}

By contrast, $\G$ often acts properly discontinuously on a contractible {\it
  complex} of smaller dimension than that of $G/K$. The minimal dimension of
such a complex, for torsion-free $\G$, is the cohomological dimension $\cd\G$
of $\G$ (with the possible exception of the case $\cd\G=2$, but such examples
are unlikely to occur among lattices in connected Lie groups) and it has been
computed in the case of arithmetic lattices by Borel and Serre
\cite{bs:corners} -- the difference $\dim G/K-\cd\G$ is equal to the $\Q$-rank
of $\G$. In many specific cases such complexes were found by Ash
\cite{ash:spine}. 

Theorem \ref{main} answers a question of Shmuel Weinberger and Kevin
Whyte. The bulk of the proof is dedicated to the special case when $G$ is a
real linear algebraic group and $\G$ is an arithmetic lattice in $G$. The
reduction to this case is presented in Section \ref{section:reduction}.
Next, we present the proof for the case when $\G=SL_n(\Z)\subset G=SL_n(\R)$,
focusing on $n=3$. This case contains most of the ideas needed in general and
has the advantage that the proof does not use the jargon of algebraic
groups. The general proof in the arithmetic case is in Section
\ref{section:arithmetic} and there are more examples preceding and following
this section. The discussion of what happens when the assumption of proper
discontinuity is replaced by the assumption that the action is isometric is in
Section \ref{section:iso}.

The method of proof is based on \cite{bkk:vk}. For $m\geq 0$ an {\it
  $m$-obstructor complex} is one that does not embed into $\R^m$ for
homological reasons. For the precise definition, see \cite{bkk:vk}. We will
  only need the fact that the join
  $$S^a*S^{k_1}_+*S^{k_2}_+*\cdots *S^{k_r}_+$$
  is an $m$-obstructor complex
  for $$m=a+k_1+k_2+\cdots+k_r+2r-1$$
and that $$S^{k_1}_+*S^{k_2}_+*\cdots *S^{k_r}_+$$ is an $n$-obstructor
complex for $$n=k_1+k_2+\cdots+k_r+2r-2$$
  where $S^k_+$ denotes the $k$-sphere
  with an extra point added, and each sphere $S^k$ is triangulated as
the join of 0-spheres. The above fact
  then follows from the Join Lemma of \cite{bkk:vk}.

When $L$ is a finite simplicial complex define the (open) cone on $L$ as
$$cone(L)=L\times [0,\infty)/L\times\{0\}.$$ If $X$ is a proper metric space,
a map $H:cone(L)\to X$ is {\it proper expanding} if it is proper and
whenever $\sigma$ and $\tau$ are disjoint simplices of $L$ then
\begin{equation}
d(H(\sigma\times [t,\infty)),H(\tau\times [t,\infty)))\to\infty
\end{equation}
as
$t\to\infty$. We also say that the cones on $\sigma$ and $\tau$ {\it
diverge} (under $H$) if (1) holds.

When $\G$ is any discrete group of type $F_\infty$ (i.e. there is an
Eilenberg-MacLane complex $K(\G,1)$ with finitely many cells in each
dimension), we write
$$\text{``}L\subset\partial\G\text{''}$$ if there is an isometric,
properly discontinuous, cocompact action of $\G$ on a proper metric
space $X$ and a proper expanding map $cone(L)\to X$. When $\G$ is a
lattice in $G$, for $X$ we will frequently use a bounded neighborhood
of a $\G$-orbit in $G$ or $G/K$ (equipped with $G$-invariant
Riemannian metrics). The formal definition of the concept ``$L\subset
\partial\Gamma$'' in \cite{bkk:vk} is more complicated and applies to
any finitely generated group, but all groups we consider in this paper
will be of type $F_\infty$ and for such groups the definition given
above is equivalent to the general definition (see \cite[Remark
11]{bkk:vk}).

The {\it obstructor
  dimension} $\obdim\Gamma$ of $\Gamma$ is defined in \cite{bkk:vk} to be
$m+2$ where $m$ is the largest integer such that ``$L\subset\partial\Gamma$''
  for some $m$-obstructor complex $L$. Passing to subgroups of finite index
  and quotients by finite normal subgroups does not change
  $\obdim$. 

The main theorem of \cite{bkk:vk} is:

\begin{thm}\cite{bkk:vk}\label{obdim}
If $\obdim\G=m+2$ then $\G$ cannot
act properly discontinuously on a contractible manifold of dimension $<m+2$.
\end{thm}

The theorem we actually prove in this paper is:

\begin{thm}\label{main2}
$\obdim\G=\dim G/K$.
\end{thm}

\begin{example} \label{uniform}
  If $\G$ is a uniform lattice in $G$ and if $G$ has finite center,
  then $\G$ acts cocompactly on the symmetric space $G/K$ which is a
  contractible manifold of dimension say $m+2$ with a complete
  Riemannian metric of nonpositive curvature (see e.g. \cite[Theorem
  V.3.1]{helgason:book}). Thus the exponential map gives a proper
  expanding map $\R^{m+2}=cone(S^{m+1})\to G/K$ (where $S^{m+1}$ is
  viewed as the unit tangent space at a point of $G/K$).  Thus
  $\text{``}S^{m+1}\subset \partial\G\text{''}$ and $\obdim\G=\dim
  G/K$, proving the theorem in this case. The same argument works for
  groups of the form $G\times A$ where $G$ is as above and $A$ is
  connected abelian.
\end{example}

Theorem \ref{main} follows from Theorem \ref{main2} and Theorem \ref{obdim}.
It is clear from Theorem \ref{obdim} that
$\obdim\G\leq\dim G/K$. To prove the reverse inequality we have to construct
an $m$-obstructor complex $L$ with $\text{``}L\subset\partial\G\text{''}$ and $m=\dim G/K-2$.

We would like to thank Dragan Mili\v ci\' c for his help with the theory of
Lie groups and algebraic groups and to Misha Kapovich and Bruce Kleiner for
several useful conversations.

\section{A fibration lemma}

A discrete version of the following lemma is in \cite{bkk:vk}. 
The proof, although it closely parallels the discrete version, is more
transparent. 

\begin{lemma}\label{ses}
  Let $G$ be a connected Lie group and $H\subset G$ a closed subgroup
  such that the homogeneous space $Q=G/H$ is contractible. The metric
  on $G$ is chosen to be invariant under left translations, on $Q$ so
  that the quotient map $\pi:G\to Q$ is a Lipschitz map, and on $H$ it
  is the subspace metric. Suppose that $H$ and $Q$ admit proper
  expanding maps of cones on finite complexes $K_H,K_Q$
  respectively. Then $G$ admits a proper expanding map of the cone on
  the join $K_H*K_Q$. Moreover, it can be arranged that the image of
  this map is contained in the $\pi$-preimage of the image of $cone(K_Q)$.
\end{lemma}

\begin{proof}
The quotient map $\pi:G\to Q$ is a fibration, and since $Q$ is contractible
there is a continuous section $s:Q\to G$.
Let $\alpha:cone(K_H)\to H$ and $\beta:cone(K_Q)\to Q$ be the given
maps. Define $$f:cone(K_H*K_Q)=cone(K_H)\times cone(K_Q)\to G$$ by
$$f(x,y)=\alpha(x)\cdot s\beta(y).$$ We have the commutative diagram
$$\begin{matrix} 
cone(K_H)&\overset{\alpha}{\longrightarrow}&H\\
\downarrow&&\downarrow\\
cone(K_H)\times cone(K_Q)&\overset{f}{\longrightarrow}&G\\
\downarrow&&\downarrow\pi\\
cone(K_Q)&\overset{\beta}{\longrightarrow}&Q\end{matrix}
$$
{\it Claim.} $f$ is a proper map.

Indeed, let $(x_i,y_i)$ be a sequence in $cone(K_H)\times cone(K_Q)$ leaving
every compact set. If the sequence $\pi f(x_i,y_i)=\beta(y_i)\in Q$ leaves
every compact set, the same is true for $f(x_i,y_i)\in G$. Otherwise, after
passing to a subsequence, we may assume that the sequence $\pi
f(x_i,y_i)=\beta(y_i)\in Q$ stays in a compact set $D\subset Q$. Then
$s\beta(y_i)$ stays in the compact set $s(D)$. Since $\beta$ is a proper map,
the sequence $y_i\in cone(K_Q)$ stays in a compact set, and thus the sequence
$x_i\in cone(K_H)$ leaves every compact set. Since $\alpha$ is a proper map,
we see that the sequence $f(x_i,y_i)=\alpha(x_i)\cdot s\beta(y_i)$ leaves
every compact set. 

{\it Claim.} If $\sigma=\sigma_H*\sigma_Q$ and $\tau=\tau_H*\tau_Q$ are
disjoint simplices of $K_H*K_Q$, then $f|cone(\sigma)$ and $f|cone(\tau)$
diverge. 

Indeed, let $(x_i,y_i)$ and $(x_i',y_i')$ be sequences in
$cone(\sigma_H)\times cone(\sigma_Q)$ and $cone(\tau_H)\times cone(\tau_Q)$
respectively, leaving every compact set. Since $\pi$ is a Lipschitz map,
if one of two sequences $\pi f(x_i,y_i)=\beta(y_i)$ and $\pi
f(x_i',y_i')=\beta(y_i')$ leaves every compact set in $Q$, then
$d_Q(\beta(y_i),\beta(y_i'))\to\infty$ (since $\beta|cone(\sigma_Q)$ and
$\beta|cone(\tau_Q)$ diverge) and consequently
$d_G(f(x_i,y_i),f(x_i',y_i'))\to\infty$. Now assume that both sequences
$\beta(y_i)$ and $\beta(y_i')$ are contained in a fixed compact set $D\subset
Q$. Then we have
\begin{equation*}
\begin{split}
d_G(f(x_i,y_i),f(x_i',y_i'))=d_G(\alpha(x_i)\cdot
s\beta(y_i),\alpha(x_i')\cdot s\beta(y_i'))=\\
d_G(1,s\beta(y_i)^{-1}\alpha(x_i)^{-1}\alpha(x_i')s\beta(y_i'))
\end{split}\end{equation*}
Since $s\beta(y_i)$ and $s\beta(y_i')$ stay in a compact set and
$$d_H(1,\alpha(x_i)^{-1}\alpha(x_i'))=d_H(\alpha(x_i),\alpha(x_i'))\to\infty$$
it follows that $$d_G(1,\alpha(x_i)^{-1}\alpha(x_i'))\to\infty$$ and
$$d_G(1,s\beta(y_i)^{-1}\alpha(x_i)^{-1}\alpha(x_i')s\beta(y_i'))\to\infty.$$ 
\end{proof}

\begin{remark}
  Another reasonable choice of a metric on $H$ would be an $H$-invariant
  Riemannian metric. Say two proper metrics $d_1$ and $d_2$ on a space $X$ are
  equivalent if there is a homeomorphism $\varphi:[0,\infty)\to [0,\infty)$
  such that $d_2(x,y)\leq \varphi(d_1(x,y))$ and $d_1(x,y)\leq
  \varphi(d_2(x,y))$ for all $x,y\in X$. Any proper expanding map $cone(L)\to
  X$ with respect to $d_1$ is also proper expanding with respect to $d_2$. For
  example, any two Riemannian $G$-invariant metrics on $G$ or on $G/K$ are
  equivalent. Various choices of metrics on $H$ as indicated above are all
  equivalent.
\end{remark}

The following consequence can be viewed as the analog of Theorem
\ref{main} in the context of nilpotent groups. Of course, the first
three statements are well known (see \cite{raghunathan:book}).

\begin{cor}\label{nilpotent}
Let $\G$ be a lattice in a simply connected nilpotent Lie group $N$. Then 
\begin{itemize}
\item $\G$
is cocompact in $N$, 
\item $N$ is diffeomorphic to Euclidean space of dimension $m+2$,
say, 
\item $N$ contains no nontrivial compact 
subgroups, and
\item $\text{``}S^{m+1}\subset\partial\G\text{''}$ and hence $\obdim\G=\dim N$.
\end{itemize}
\end{cor}

\begin{proof}
Let $Z$ be the center of $N$. By \cite[Proposition
  2.17]{raghunathan:book} the intersection $Z\cap \G$ is a lattice in $Z$. All
  claims now follow by induction on $\dim N$ from the exact sequence $$1\to
  Z\to N\to N/Z\to 1$$ and Lemma \ref{ses}.
\end{proof}

\begin{example}\label{heisenberg}
  Let $N_n$ be the group of real upper-triangular matrices with 1's on the
  diagonal, let $N_n(\Z)$ be the lattice in $N_n$ consisting of matrices with
  integral entries, and let $m+2=1+2+\cdots +(n-1)$ denote the number
  of matrix positions
  above the diagonal. Then $$\text{``}S^{m+1}=*_{i=1}^{m+2}S^0\subset\partial
  N_n(\Z)\text{''}.$$
  Specifically, if we regard the cone on $S^{m+1}$ as $\R^{m+2}$
  with a coordinate for every position above the diagonal, and we regard the cone
  on a simplex of $S^{m+1}$ as the set of matrices in $H_n$ where the entries
  in the specified positions have specified signs and the other entries above the
  diagonal are 0, then this map $cone(S^{m+1})\to H_n$ is proper and
  expanding.
\end{example}

Another application of Lemma \ref{ses} is to semi-simple groups with
infinite center, e.g. $\widetilde {SL_2(\R)}$.

\begin{cor}\label{infinite center}
  Suppose $G$ is a connected semi-simple Lie group with infinite
  center $Z$, and let $\G\subset G$ be a lattice such that $\G\cap Z$
  has finite index in $Z$ so that $\G/(\G\cap Z)$ is a lattice in
  $G/Z$. Then $\obdim \G\geq \obdim (\G/\G\cap Z)+\rank(Z)$.
\end{cor}

\begin{proof} Let $K$ be a maximal compact subgroup in $G/Z$. Now apply Lemma
  \ref{ses} to the preimage $H\subset G$ of $K$ in $G$.
\end{proof}

\section{Reduction to arithmetic lattices}\label{section:reduction}

We will review the terminology and basic facts about algebraic and
arithmetic groups in Section \ref{section:arithmetic}.
The bulk of the paper is dedicated to the proof of the following theorem:

\begin{thm} \label{arithmetic}
  Let $G\subset GL_n(\C)$ be a semi-simple linear algebraic group defined over
  $\Q$ and let $\G_\Z=G\cap GL_n(\Z)\subset G_\R$ be the standard arithmetic
  lattice in the group of real points.  Then $\obdim \G_\Z=\dim G_\R/K$.
\end{thm}

In this section we deduce Theorem \ref{main2} from Theorem \ref{arithmetic}. We
first prove another special case.

\begin{thm} \label{marg}
Let $G$ be a semi-simple real algebraic group and $\G$ a lattice in $G$. Then
$\obdim \G_\Z=\dim G_\R/K$. 
\end{thm}

\begin{proof}
  We can assume that the component of the identity of $G_\R$ has no
  compact factors. Further, using the Product Lemma \cite{bkk:vk}, we
  may assume that $\G$ is an irreducible lattice in $G_\R$. If the
  real rank of $G$ is $>1$ then the celebrated theorem of Margulis
  \cite[Theorem 6.1.2]{zimmer:book} says that $\G$ is an arithmetic
  lattice (with respect to some $\Q$-structure on $G$) and the theorem
  follows from Theorem \ref{arithmetic}. Now suppose that the real
  rank of $G$ is 1. Then the symmetric space $G/K$ is real, complex,
  or quaternionic hyperbolic space or the Cayley plane (see
  \cite[Chapter X]{helgason:book}). If $\G$ is a uniform lattice
  acting cocompactly on $G/K$, then the statement follows from Example
  \ref{uniform}, so we may assume that $\G$ is a nonuniform
  lattice. If $G/K$ is the real hyperbolic space $\H^{m+2}$ let $P$ be
  a maximal parabolic subgroup of $\G$. Then $P$ is commensurable to
  $\Z^{m+1}$ and there is a proper expanding map of the cone on $S^m$
  into a horosphere stabilized by $P$. Adding a ray that diverges from
  this horosphere but stays within a bounded distance from a
  $\G$-orbit produces a proper expanding map from the cone on $S^m_+$
  and shows that $\obdim(\G)=m+2$. If $G/K$ is the complex hyperbolic
  space of complex dimension $d=\frac{m+2}{2}$, then a maximal
  parabolic subgroup $P$ of $\G$ is a lattice in the Heisenberg group
  $N$ that is a central extension $1\to\R\to N\to \C^{d-1}\to 1$. By
  Corollary \ref{nilpotent} there is a proper expanding map of the
  cone on $S^m$ into a horosphere stabilized by $P$. The argument now
  follows as in the real case. If $G/K$ is quaternionic hyperbolic
  space or the Cayley plane, then $\G$ is an arithmetic lattice
  \cite{gs:harmonic} and the claim again follows from Theorem
  \ref{arithmetic}.
\end{proof}

\begin{proof}[Proof of Theorem \ref{main2}]
  First suppose that the center $Z$ of $G$ is finite. If necessary, replace
  $G$ by $G/Z$ so that the center is trivial. Then $G$ has the structure of
  (the identity component of) a
  linear real algebraic group (see e.g. \cite[Proposition 3.1.6]{zimmer:book})
  and the proof is reduced to Theorem \ref{marg}.
  
  Now suppose that the center $Z$ of $G$ is infinite. Again we may assume that
  $G$ has no compact factors. Then by \cite[5.17]{raghunathan:book} $Z\Gamma$
  is a discrete subgroup of $G$, so \cite[1.13]{raghunathan:book} (with $H=Z$)
  implies that $\G\cap Z$ has finite index in $Z$. It now follows from
  Corollary \ref{infinite center} and the centerless case applied to the
  lattice $\G/(Z\cap\Gamma)\subset G/Z$ that
  $$\obdim\Gamma\geq\rank Z+\obdim \G/(Z\cap\Gamma)=\rank Z+\dim (G/Z)/C$$
  where $C$ is a maximal compact subgroup of $G/Z$. It remains to show that
  $$\dim G/K=\rank Z+\dim (G/Z)/C.$$ After passing to finite covers,
we may assume that $C=T\times C'$ where $T$ is a torus and $C'$ is a
simply connected compact group (see e.g. \cite[Theorem
4.29]{tonci}). The preimage of $C$ in $G$ decomposes as $\R^m\times
T'\times C'$ where $T'$ is a torus and $m=\rank Z$. 
Then $K$ can be identified with $T'\times C'$, so we have
\begin{equation*}
\begin{split}
\dim G/K=\dim G-\dim K=\dim (G/Z)-(\dim C'+\dim T')=\\
\dim
(G/Z)-\dim C+m=\dim (G/Z)-\dim C+\rank Z\end{split}\end{equation*}
\end{proof}

\section{Isometric Actions}\label{section:iso}

We need the following fact from the theory of Lie groups.

\begin{lemma}\label{rep}
If $H$ is a connected noncompact Lie group, then there is a
representation $\rho:H\to GL_N(\R)$ such that the closure of the image
$\overline{\rho(H)}$ is noncompact.
\end{lemma}

\begin{proof}
Assume that the adjoint representation has precompact image. Then the
Lie algebra $\frak h$ of $H$ admits an $ad$-invariant inner product and
hence breaks up as the direct sum of simple and abelian Lie
algebras. The simple summands have compact type, so the integral
subgroup $C\subset H$ corresponding to the sum of all simple summands
is compact, as well as normal in $H$. The quotient $H/C$ is a
noncompact connected abelian group and it therefore maps onto $\R$.
\end{proof}

Note that the group $GL_N(\R)$ above can be replaced by $PGL_{N+1}(\C)$
since there is a proper map $GL_N(\R)\to PGL_{N+1}(\C)$.

\begin{thm}
  Suppose that $M$ is a contractible manifold equipped with a
  complete Riemannian metric. Let $G$ be a connected simple Lie group
  of rank $>1$ with finite center and let
  $\G$ be a lattice in $G$. If $\dim M<\dim G/K$ then
  every isometric action of $\G$ on $M$ has a bounded orbit. In
  particular, if the metric on $M$ is CAT(0), then every isometric
  action of $\G$ on $M$ has a global fixed point.
\end{thm}

\begin{proof}
  The isometry group $Isom(M)$ of $M$ is a
  Lie group.  Let $\varphi:\G\to Isom(M)$ be the given action, and
  denote by $H$ the closure of $\varphi(\G)$. Then $H$ is also a Lie
  group. If $H$ has infinitely many components, then the
  Margulis-Kazhdan theorem \cite[Theorem 8.1.2]{zimmer:book} implies
  that $\varphi:\G\to Isom(M)$ has finite kernel and the image of
  $\varphi$ is closed. It follows that the action of $\G$ on $M$ is
  properly discontinuous, contradicting the assumption $\dim M<\dim
  G/K$ and Theorem \ref{main}. Now suppose that $H$ has only finitely
  many components. After passing to a subgroup of $\G$ of finite
  index, we may assume that $H$ is connected. If $H$ is compact, then
  the $H$-orbits, and hence $\G$-orbits, are bounded. Thus assume that
  $H$ is noncompact. Then $H$ admits a representation $\rho:H\to
  PGL_{N+1}(\C)$ whose image has noncompact closure (see Lemma
  \ref{rep}). We can arrange in addition that $\rho$ is trivial on the
  finite central subgroup $\varphi(\G\cap Z(G))$ (by applying Lemma
  \ref{rep} to $H/\varphi(\G\cap Z(G))$). By the Margulis'
  super-rigidity \cite[Theorem 5.1.2]{zimmer:book} (applied to an
  algebraic $\R$-group that has $G/Z(G)$ as the identity component
  \cite[Proposition 3.1.6]{zimmer:book} and to the lattice $\G/(\G\cap
  Z(G))$) there is a continuous extension $\tilde \rho:G\to
  PGL_{N+1}(\C)$. The kernel $Ker(\tilde \rho)\neq G$ is a normal
  subgroup of $G$ and is therefore contained in the center. It follows
  from \cite{goto} (see also \cite[Exercise II.D.4]{helgason:book})
  that $\tilde \rho(G)$ is a closed subgroup of $PGL_{N+1}(\C)$, thus
  the image of $\varphi$ is closed and discrete and $H$ cannot be
  connected, contradiction. The last sentence follows from the Cartan
fixed point theorem
  \cite[Corollary II.2.8]{bh:book}.
\end{proof}

If $G$ has rank 1, one can say the following. When $G/K$ is real or complex
hyperbolic space there exist lattices with positive first Betti number, and
they of course act isometrically on $\R$ with unbounded orbits. In the
real case this is a result of Millson \cite{millson:betti} and the complex
case a result of Kazhdan (see \cite[Corollary VIII 5.9]{bw:book}). The
cases of quaternionic hyperbolic space or the Cayley plane are open. By
Property T, lattices in these cases have trivial first Betti number. One
possible scenario would be to have infinite quotients of smaller cohomological
dimension that act discretely on a contractible manifold of smaller
dimension. 
\vskip .2cm
\noindent
{\bf Question:} Do lattices in $Sp(n,1)$ have infinite quotients of strictly
smaller (virtual or rational) cohomological dimension?
\vskip .2cm
If $G$ is not simple, but only semi-simple (and with finite center),
then $G/K$ decomposes (see \cite[Proposition V.4.2, Proposition
VIII.5.5]{helgason:book}) as the direct product $X_1\times X_2\times
\cdots\times X_m$ of irreducible symmetric spaces and an irreducible
lattice $\G\subset G$ acts isometrically on each $X_i$ with unbounded
orbits. After replacing (without loss of generality!) $G$ by a finite
cover, $G$ also decomposes as $G_1\times G_2\times\cdots\times
G_m\times C$ where $C$ is a compact group, $G_i$ are simple and
$G_i/K_i=X_i$. 

\begin{thm}
Suppose that $G=G_1\times G_2\times\cdots\times G_m$ ($m>1$) is the
product of noncompact simple Lie groups with finite center and that an
irreducible lattice $\G\subset G$ acts isometrically on a piecewise
Riemannian contractible manifold $M$ with $\dim M<\dim X_i$ for all
$i=1,2,\cdots,m$. Then all $\G$-orbits are bounded.
\end{thm} 

\begin{proof}
Follow the
proof above {\it verbatim} until the extension $\tilde\rho:G\to
PGL_{N+1}(\C)$ is considered. Now the kernel is a closed normal
subgroup and its Lie subalgebra can be assumed, after reordering the
factors, to be equal to the Lie subalgebra of $G_1\times\cdots\times
G_i$ for some $1\leq i<m$ (if the kernel is discrete the argument
concludes the same way). It follows that the product 
$$G\to G_1\times\cdots\times G_i\times PGL_{N+1}(\C)$$ of projections
and of $\tilde\rho$ has finite kernel. Now \cite{goto} again implies that the
image of this map is closed, and hence the diagonal action of $\G$ on
$X_1\times X_2\times\cdots\times X_i\times M$ is properly
discontinuous. Theorem \ref{main} now implies that $\dim M\geq \dim
X_{i+1}+\cdots+ \dim X_m\geq \dim X_m$.
\end{proof}

\begin{remark}
The metric on $M$ can be allowed to be more general than Riemannian in
both theorems above. The proof only requires that the metric is proper
and that $Isom(M)$ is a Lie group. The latter is satisfied e.g. if the
Hausdorff dimension of the metric is $<\dim M+2$ by
\cite{repovs:hilbert}. Hilbert-Smith conjecture implies that $Isom(M)$
is a Lie group for any metric on $M$.
\end{remark}

\section{Examples}

In the examples that follow (before and after the proof of Theorem
\ref{arithmetic}) $G$ will be a semi-simple group of matrices, $\G$ a
lattice in $G$ and $K$ a maximal compact subgroup of $G$. If the
dimension of the symmetric space $G/K$ is $m+2$, the goal is to
construct a proper expanding map $cone(L)\to G$ of an $m$-obstructor
complex $\L$ in $G$ whose image is in a bounded neighborhood of
$\G$. To compute the dimension of $X=G/K$ we will make use of the
Iwasawa decomposition $G=KAN$ that implies $\dim X=\dim A+\dim N$.

\section{$SL_3(\Z)$}

We now discuss the case $\G=SL_3(\Z)$ in detail.  Here $G=SL_3(\R)$,
$K=SO_3$, and the
symmetric space $X=G/K$ is 5-dimensional.
We will construct the 3-obstructor complex $\L=S^0_+*S^1_+$ in
$\partial\G$.

\subsection{The cuspidal complex $C$}
We will focus on the subgroup $B\subset \G$ consisting of upper-triangular
matrices with 1's on the diagonal and on the 6 conjugates of $B$ obtained by
simultaneously permuting rows and columns. As we know (see Example
\ref{heisenberg}) each copy of $B$ has a natural 2-sphere at infinity. To
understand how different 2-spheres fit together we define the simplicial
complex $C$. Its vertices are the 6 off-diagonal positions in a
$3\times 3$ matrix. A vertex $ij$ can be viewed as an oriented arrow from a
point labeled $i$ to a point labeled $j$ ($i,j\in\{1,2,3\}$). A collection of
vertices of $C$ spans a simplex if there is a conjugate of $B$ as above that
has non-zero entries in all the corresponding positions. In the language of
arrows, the condition is that there are no oriented cycles. See Fig. 1.
\begin{figure}[htbp]
\begin{center}

\input{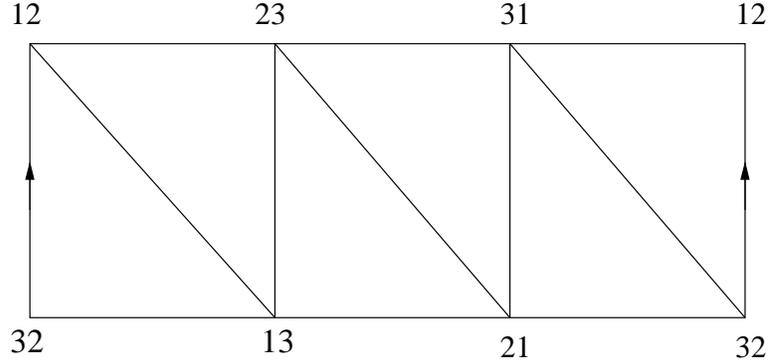}

\caption{\bf The complex $C$ is an annulus with a drum triangulation}
\label{figure:c}
\end{center}
\end{figure}

\subsection{The complex $SC$}
We now define a functorial procedure that assigns a complex $SC$ to
the complex $C$. Each simplex will be replaced by the sphere of the
same dimension. A vertex $v$ of $C$ corresponds to two vertices,
thought of as $v+$ and $v-$, in $SC$. Thus the vertex set $SC^{(0)}$
of $SC$ is the vertex set $C^{(0)}$ of $C$ crossed with
$\{+,-\}$. There is a natural projection $\pi:SC^{(0)}\to C^{(0)}$
that forgets the sign. A collection of vertices of $SC$ spans a
simplex iff $\pi$ is injective on this collection and the image in
$C^{(0)}$ is the vertex set of a simplex. The projection $\pi$ extends
to a simplicial map $\pi:SC\to C$. A simplex of dimension $k$ in $C$
has precisely $2^{k+1}$ lifts to $SC$ and the full preimage is a
$k$-sphere triangulated as the $(k+1)$-fold join of 0-spheres.

\subsection{Two lemmas}

\begin{lemma}[Lemma A]
$SC$ contains a 3-obstructor complex $\L$.
\end{lemma}

\begin{lemma}[Lemma B]
$\text{``}\L\subset \partial SL_3(\Z)\text{''}$
\end{lemma}

\begin{remark}
We don't
know if $\text{``}SC\subset \partial SL_n(\Z)\text{''}$.
\end{remark}

We will prove Lemma B for $SL_n(\Z)$ in the next section. For Lemma A,
consider the subcomplex $\L$ of $SC$ formed by the full preimage (which is
a 2-sphere) of one of
the triangles in $C$, say $<12,23,13>$; add to it an equatorial disk
which is half of the sphere $\pi^{-1}(<13,23,21>)$; finally add the
vertex $32+$ and edges joining it to the north and south poles and the
center of the equatorial disk. See Fig. 2.

\begin{figure}[htbp]
\begin{center}

\input{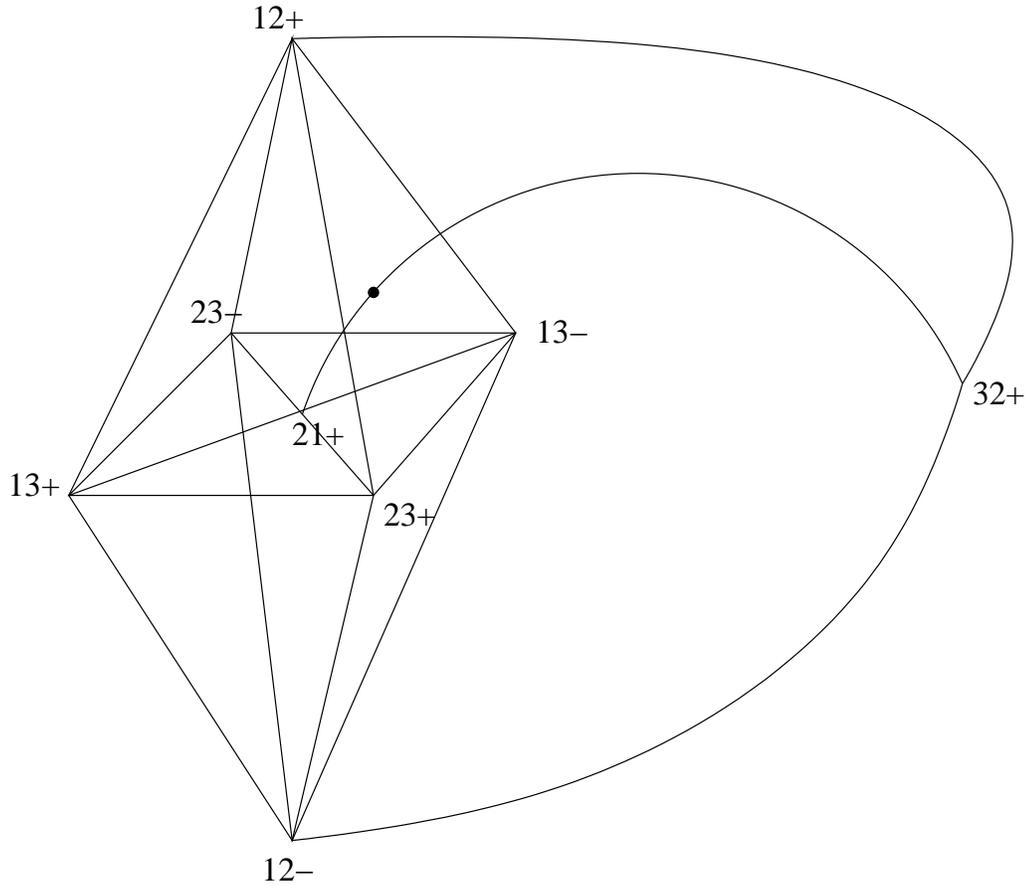}

\caption{\bf A 3-obstructor subcomplex $\L\subset SC$ nearly embedded in
  $\R^3$. The Van Kampen
  intersection point is indicated.}
\label{figure:sc}
\end{center}
\end{figure}

This subcomplex $\L$ is the join of the disjoint union of the circle
$\pi^{-1}(<13,23>)$ and the vertex $32+$ and the three points $12+,12-,21+$.
By the Join Lemma of \cite{bkk:vk}, $\L$ is a 3-obstructor complex. This shows
$\obdim SL_3(\Z)=5$.

\section{$\G=SL_n(\Z)$}

Here $G=SL_n(\R)$, $A$ is the group of diagonal matrices in $G$, and
$N$ is the group of upper triangular matrices in $G$ with 1's on the
diagonal. Thus $\dim X=\dim
AN=(1+2+\cdots+n)-1=\frac{n^2}{2}+\frac{n}{2}-1$. We let
$$\L=S^0_+*S^1_+*\cdots *S^{n-2}_+$$ which is an $m$-obstructor
complex for $m=0+1+\cdots+(n-2)+2(n-2)=\frac{n^2}{2}+\frac{n}{2}-3$.

We again define the complex $C$ with $n^2-n$ vertices, one for every
off-diagonal position of an $n\times n$ matrix, and a collection of
vertices spans a simplex iff the corresponding positions are all above the
diagonal after a simultaneous permutation of rows and columns. Let
$SC$ be the corresponding complex whose vertices are signed vertices
of $C$ and whose simplices are lifts of simplices of $C$.

\begin{lemma}[Lemma A]
$SC$ contains the $m$-obstructor complex $\L$.
\end{lemma}

\begin{lemma}[Lemma B]
$\text{``}\L\subset \partial SL_n(\Z)\text{''}$
\end{lemma}

\begin{proof}[Proof of Lemma A]
The sphere $S^{k-2}$ is the full preimage of the $(k-2)$-simplex
$<1k,2k,\cdots,(k-1)k>$ of $C$ and we add the vertex $k(k-1)+$ to form
$S^k_+$. It is straightforward to check that the subcomplex of $SC$
spanned by the described vertices is precisely a copy of $\L$.
\end{proof}

\subsection{Lemma B -- discussion}
We first define a map from the cone on $SC$ into $G=SL_n(\R)$. 
A simplex of $SC$ corresponds to a collection of positions with signs and the
first guess might be that the cone on such a simplex is sent to the
subset of $G$ having 1's on the diagonal, entries of appropriate sign
in the positions corresponding to the vertices of the simplex, and 0's in
all other positions. The problem with this rule is that cones on disjoint
simplices don't diverge. For example (letting $R$ be large), let
$$
A=\left[
{\begin{array}{rcc}
1 & R & 1 \\
0 & 1 & 1 \\
0 & 0 & 1
\end{array}}
 \right]  
$$ and
$$
B=  \left[
{\begin{array}{rcc}
1 & 0 & -R \\
0 & 1 & 0 \\
0 & 1 & 1
\end{array}}
 \right]  
$$ Then $$A^{-1}B=\left[ {\begin{array}{rcc} 
1 & -1 & -1 \\ 
0 & 0 & -1\\ 
0 & 1 & 1
\end{array}}
 \right]  
$$ so $A$ and $B$ are bounded distance apart. Instead, we separate the
positions below diagonal from the positions above the diagonal, thus obtaining
two matrices, one upper-triangular and the other lower-triangular. We
then send the point in the cone to the product of the two
matrices. For example, the cone on $<23+,31+>$ is thought of not as
the set of matrices
$$\left[
{\begin{array}{rcc}
1 & 0 & 0 \\
0 & 1 & x \\
y & 0 & 1
\end{array}}\right]$$ with $x,y\geq 0$ but as the set of matrices
$$\left[
{\begin{array}{rcc}
1 & 0 & 0 \\
0 & 1 & x \\
0 & 0 & 1
\end{array}}\right]
\left[
{\begin{array}{rcc}
1 & 0 & 0 \\
0 & 1 & 0 \\
y & 0 & 1
\end{array}}\right]=
\left[
{\begin{array}{rcc}
1 & 0 & 0 \\
xy & 1 & x \\
y & 0 & 1
\end{array}}\right]$$ with $x,y\geq 0$.

In
the example above, $A$ is already upper-triangular and would represent
the image of the chosen point in the cone on $<12+,23+,13+>$, but the
matrix $B$ would be replaced by the product 
$$\left[
{\begin{array}{rcc}
1 & 0 & -R \\
0 & 1 & 0 \\
0 & 0 & 1
\end{array}}\right]\left[
{\begin{array}{rcc}
1 & 0 & 0 \\
0 & 1 & 0 \\
0 & 1 & 1
\end{array}}\right]=\left[
{\begin{array}{rcc}
1 & -R & -R \\
0 & 1 & 0 \\
0 & 1 & 1
\end{array}}\right]
$$
If this matrix is denoted by $B'$, then we compute that $$A^{-1}B'=
\left[
{\begin{array}{rcc}
1 & -R-1 & -1 \\
0 & 0 & -1 \\
0 & 1 & 1
\end{array}}\right]
$$
so $A$ and $B'$ are far apart.  In this case it can be shown that this
first guess for the map suffices to prove that $\text{``}\L\subset
\partial SL_n(\Z)\text{''}$; however, in general we find it easier to
deal with the map corresponding to separating the upper-triangular
from the lower-triangular part.

\begin{lemma} \label{2.5}
  For every $n$ there is a function $\phi:[0,\infty)\to [0,\infty)$ with the
  following properties.  Let $\Lambda=(l_{ij}),\Lambda'=(l'_{ij})$ be two
  lower triangular matrices with 1's on the diagonal, and let $U,U'$ be two
  upper-triangular matrices with 1's on the diagonal. Assume:
\begin{itemize}
\item If $l_{ij}\neq 0$ or if $l'_{ij}\neq 0$ then either $i=j$ or
$i=j+1$.
\item $l_{j+1,j}l'_{j+1,j}=0$ for all $j$.
\end{itemize}
If at least one of $\Lambda,\Lambda'$ contains an entry of absolute value
$>\phi(R)$ then $S=(U\Lambda)^{-1}U'\Lambda'$ has an entry of absolute value
$>R$.
\end{lemma}

\begin{proof}
  We can write $U\Lambda=U'\Lambda'S$. Assume that the statement is
  false, and let $U_k,U_k',\Lambda_k,\Lambda_k',S_k$ be a sequence of
  examples with all entries of $S$ bounded by $R$ and some entries of
  $\Lambda$ or $\Lambda'$ going to infinity. By passing to a
  subsequence we may assume that all entries converge in
  $[-\infty,\infty]$. Call a position of a matrix `small' if the
  corresponding entries stay bounded and otherwise call it
  `large'. For example, all positions of $S$ are small. Note that by
  Cramer's Rule ($\det S=1$) all positions of $S^{-1}$ are also small. For
  simplicity, we drop the subscripts from $\Lambda_k,U_k,\cdots$.

Consider the lowest nonzero off-diagonal entry $x$ of $\Lambda$ or
$\Lambda'$, i.e. $l_{j+1,j}$ or $l'_{j+1,j}$ such that
$l_{m+1,m}=l'_{m+1,m}=0$ for all $m>j$. There are two cases:

{\bf Case 1:} The position of $x$ is small. We can replace $x$ by 0
without affecting other entries of the same matrix by multiplying on
the right by the elementary lower-triangular matrix that has entry
$-x$ in the position $j+1,j$. We can compensate this operation by either
multiplying $S$ on the right by the same matrix (in case that the
entry $x$ belongs to $\Lambda$) or by multiplying $S$ on the left by the
inverse of this matrix (in case $x$ belongs to $\Lambda'$). All entries of
the new matrix $S$ are still small (though not bounded by $R$
perhaps).  Now proceed to the next lowermost off-diagonal entry of $\Lambda$
or $\Lambda'$.

{\bf Case 2:} The position of $x$ is large.
Say $x$ belongs to $\Lambda$ (the other case is analogous). We write 
$${U'}^{-1}U\Lambda=\Lambda'S$$
The right hand side is obtained from the matrix $S$ with small
entries 
by applying elementary row operations in which a multiple of the
$p^{th}$ row is added to the row below for certain $p< j$. Therefore
all rows bellow row $j$ of the matrix $\Lambda'S$ are equal to the
corresponding entries of the small matrix $S$.

The left hand side is obtained from $\Lambda$ by applying elementary row
operations in which a multiple of a row is added to some row
above. It follows that the entry of ${U'}^{-1}U\Lambda$ in position $(j+1,j)$ is
$x$, which is large. Contradiction
\end{proof}

\begin{proof}[Proof of Lemma B]
  Let $U\Lambda$ and $U'\Lambda'$ be two matrices representing points in cones
  on disjoint simplices of $\L$. Assuming that the two points are far away
  from the cone point, we have to show that the distance between $U\Lambda$
  and $U'\Lambda'$ is large. If one of $\Lambda,\Lambda'$ has a large entry,
  this follows from Lemma \ref{2.5}. If both $\Lambda$ and $\Lambda'$ have
  small entries, then $U\Lambda$ is close to $U$, while $U'\Lambda'$ is close
  to $U'$, so we must show that $U$ and $U'$ are far apart. In this case $U$
  and $U'$ will have large entries, so the statement follows from Example
  \ref{heisenberg}.
\end{proof}

We can see more concretely how the cone on $L$ is mapped to $G$.
The sphere $S^{p-2}$ is
identified with the unit sphere in $\R^{p-1}$, which in turn is mapped to $G$
via $$(x_1,x_2,\cdots,x_{p-1})\mapsto \left[
{\begin{array}{rcccccc}
1&&&&x_1&& \\
&1&&&x_2&& \\
&&&&&& \\
&&&&x_{p-1}&& \\
&&&&1&& \\
&&&&&& \\
&&&&&&1\end{array}}
 \right]$$
and the extra point determines the ray of positive numbers in the $(p,p-1)$-position.

\section{$\G=SL_n(\Z[\sqrt{2}])$}

$\O=\Z[\sqrt{2}])$ is a lattice in $\R^2$ under the embedding
$a+b\sqrt{2}\mapsto (a+b\sqrt{2},a-b\sqrt{2})$ and likewise $\G$ is a
lattice in $G=SL_n(\R)\times SL_n(\R)$ (with the same map applied to
each entry). Note that $\O$ is a ring and the group of units $\O^*$
has rank 1 (e.g. $1+\sqrt{2}$ has infinite order).  The symmetric
space is $X\times X$ with $X=SL_n(\R)/SO_n$. We define $\L$ as the
join
$$\L=S^{n-2}*S^1_+*S^3_+*\cdots *S^{2n-3}_+.$$
Then $\L$ is an $m$-obstructor complex
for $m=(n-1)+1+3+\cdots+(2n-3)+2(n-2)=n^2+n-4$ while $X\times X$ has dimension
$n^2+n-2$.

We define the map $\Psi:cone(\L)\to G$ analogously to the $SL_n(\Z)$
construction, with a new feature that there are diagonal entries
different from 1 this time.
 
Fix some $p\in\{2,3,\cdots,n\}$ and consider the subgroup $\G_p$ of $\G$
consisting of matrices that have any entry from $\Z[\sqrt{2}]$ in positions
$(1,p),(2,p),\cdots,(p-1,p)$, 1's on the diagonal, and 0's in the remaining
positions.
$$
\left[
{\begin{array}{rcccccc}
1 & & & & x_1 & & \\
&1&&&x_2&& \\
&&&&&& \\
&&&&x_{p-1}&& \\
&&&&1&& \\
&&&&&& \\
&&&&&&1
\end{array}}
 \right]  
 $$

 To this subgroup (isomorphic to $\Z^{2p-2}$) we associate the subgroup $G_p$
 of $G$ isomorphic to $\R^{2p-2}$ consisting of the set of pairs $(A,B)$ of
 matrices as above with real entries. Note that $\G_p$ is a cocompact lattice
 in $G_p$. We view $G_p$ as the cone on $S^{2p-3}$ with the cone on a simplex
 of $S^{2p-3}$ corresponding to the subset of $G_p$ where certain entries of
 $A$ and $B$ above the diagonal and in column $p$ are
 nonnegative, certain others are nonpositive, and the remaining entries are 0. 
 
 Next, we define a ray $R_p$ in $G$ to consist of pairs of matrices
 $(\Lambda,\Lambda)$ with $\Lambda$ having 0's above the diagonal, 1's on the
 diagonal, entry $x\geq 0$ in position $(p,p-1)$ and 0's in all other positions. This
 ray is also within a bounded distance from $\G$. We identify the cone on
 $S^{2p-2}_+$ with ${G_p}_+:=G_p\cup R_p$.
 
 Finally, we define $\G_0$ to consist of diagonal matrices in $\G$, and we
 define $G_0$ to consist of pairs $(D,D^{-1})\in G$ where $D$ is a diagonal
 matrix with positive diagonal entries. We have that $\G_0$ is an abelian
 group of rank $n-1$, $G_0\cong \R^{n-1}$, and $G_0$ is within a bounded
 neighborhood of $\G_0$. We identify $G_0$ with the cone on $S^{n-2}$ where
 the cone on a simplex of $S^{n-2}$ corresponds to the pairs $(D,D^{-1})$ with
 certain diagonal entries of $D$ bounded below by 1, certain others bounded
 above by 1, and the remaining diagonal entries equal to 1.

We now define the map $\Psi:cone(\L)\to G$. We are identifying $cone(\L)$ with
$G_0\times {G_2}_+\times \cdots\times {G_n}_+$. Let
$((D,D^{-1}),(A_2,B_2),\cdots,(A_{n},B_{n}))\in cone(\L)$. From this
data we first form 3 pairs of matrices: 
\begin{itemize}
\item diagonal pair $(D,D^{-1})$,
\item lower-triangular pair $(\Lambda,\Lambda)$ -- it is formed by
  superimposing (equivalently, adding entries below the diagonal) all
  lower-triangular pairs $(\Lambda_p,\Lambda_p)$ appearing in the sequence
  $(A_p,B_p)$, and
\item upper-triangular pair $(U,V)$ -- it is formed by superimposing
  (equivalently, adding entries above the diagonal) all upper-triangular pairs
  $(A_p,B_p)$ appearing in the above sequence.
\end{itemize}
We define the image under $\Psi$ of the given sequence to be the pair
$$(UD\Lambda,VD^{-1}\Lambda)\in G.$$

\begin{lemma}[Lemma B]
$\text{``}\L\subset\partial SL_n(\Z[\sqrt{2}])\text{''}$.
\end{lemma}

The proof is similar to the proof in the case of $SL_n(\Z)$ except for the
added complication of diagonal matrices. It will not be given here and
we appeal to the general case of Lemma B given is Section \ref{sec:B}.

\section{$\G=SL_n(\O)$}

Here $\O$ is the ring of integers in a number field $k$. Let $r$ be
the number of real places of $k$ and $s$ the number of
complex-conjugate places. Then, as an abelian group, $\O\simeq
\Z^{r+2s}$. The group $\O^*$ of units is a finitely generated group of
rank $r+s-1$. For example, when $\O=\Z[i]$ then $r=0$ and $s=1$, and
when $\O=\Z[\sqrt{2}]$ then $r=1$ and $s=0$. $\O$ has $r$ embeddings
in $\R$ and $s$ embeddings in $\C$, and there is an induced diagonal
embedding $\G\subset G=SL_n(\R)^r\times SL_n(\C)^s$ and $\G$ is a
lattice there. The symmetric space associated to $SL_n(\C)$ has
dimension equal to the dimension of $A$ ($=n-1$) plus the dimension of
$N$ ($=2(1+2+\cdots+(n-1))$), so it equals $n^2-1$. Thus the dimension
of the symmetric space of $G$ is
$$r(\frac{n^2}{2}+\frac{n}{2}-1)+s(n^2-1)$$

We set
$$
\begin{matrix}
\L=&S^{(r+2s)+(r+s-1)-1}_+*&S^{2(r+2s)+(r+s-1)-1}_+*S^{3(r+2s)+(r+s-1)-1}_+*\cdots
*\\ &S^{(n-1)(r+2s)+(r+s-1)-1}_+& \end{matrix}$$ which is an
$m$-obstructor complex for
$$\begin{matrix} m&=(1+2+\cdots+(n-1))(r+2s)+(n-1)(r+s-2)+2(n-2)\\
&=r(\frac{n^2}{2}+\frac{n}{2}-1)+s(n^2-1)-2\hfill \end{matrix}$$ The
complex $S^{(p-1)(r+2s)+(r+s-1)-1}_+$ is identified with the unit
sphere in the space $\R^{(p-1)(r+2s)+(r+s-1)}$ union a ray. The
Euclidean space is realized (as an abelian group with that rank) by
upper triangular matrices with entries $(1,p),(2,p),\cdots,(p-1,p)$
from $\O$, entry $(p,p)$ from $\O^*$, entry $(1,1)$ the inverse of
entry $(p,p)$, other diagonal entries are 1, and all other entries are
0. The ray is realized by matrices with diagonal entries 1, and
positive real entry in position $(p,p-1)$.

\section{Proof of Theorem \ref{arithmetic}} \label{section:arithmetic}

In this section we consider the case of arithmetic lattices. We use
\cite{raghunathan:book}, \cite{borel:lag}, and \cite{borel:ga} as
general references on algebraic and arithmetic groups. Let $G\subset
GL_N(\C)$ denote a connected semi-simple linear algebraic group defined
over $\Q$. If $k$ is a subring of $\C$ (usually $\Z$, $\Q$, $\R$, or
$\C$) denote by $G_{k}$ the group $G\cap GL_N({k})$ of $k$-points. It
is a theorem of Borel-Harish-Chandra \cite{bhc:arithmetic} that $G_\Z$
is a lattice in $G_\R$.  This is the standard arithmetic lattice. Let
$S$ be the maximal torus in $G$ which is split over $\Q$ and let $M$ be
the largest connected subgroup of the
centralizer $Z(S)$ which
is anisotropic, i.e. does not have any nontrivial $\Q$-split tori (in other
words, the $\Q$-rank of $M$ is 0).
We also have that the component of the identity
$Z(S)^0=S\cdot M$, i.e. the multiplication map $S\times M\to Z(S)^0$
is surjective and has finite kernel. Further, $M$ is reductive
\cite[IV.13.17 Corollary 2]{borel:lag}, i.e. after a finite cover it
decomposes as the product of a torus and a semi-simple $\Q$-group.

There is the usual decomposition of the Lie algebra $\frak g$ of $G$:
$$\frak g=\frak g_0\oplus \bigoplus_{\alpha\in\Phi}\g_\alpha$$ into
root spaces. Each $\alpha\in\Phi$ is a rational character $\alpha:S\to
GL_1(\C)$ and $\g_\alpha$ is the associated root space. It is
customary to use additive notation in the group of characters. Choose
an ordering on $\Phi$ and denote by $\Phi^+$ the set of positive roots
and by $\Delta$ the set of simple roots (those roots not expressible
as sums of other roots in $\Phi^+$). The cardinality $r$ of $\Delta$
is equal to the dimension of $S$ and is called the {\it $\Q$-rank} of
$G$. Lattice $G_\Z$ is cocompact in $G_\R$ if and only if the
$\Q$-rank is 0 (this statement holds even for reductive
groups). Unlike in the case of algebraically closed fields, the root
spaces $\g_\alpha$ may have dimension $>1$ and the root system $\Phi$
may not be reduced (we may have $0\neq\alpha,2\alpha\in\Phi$). Of
course, $\Phi$ might not be irreducible (the Dynkin diagram could be
disconnected). Every irreducible component of $\Phi$ is either reduced
(i.e. it is of type $A_n$ ($n\geq 1$), $B_n$ ($n\geq 3$), $C_n$
($n\geq 2$), $D_n$ ($n\geq 4$), $E_6$, $E_7$, $E_8$, $F_4$, or $G_2$)
or unreduced (i.e. it is of type $BC_n$ ($n\geq 1$)). See
e.g. \cite{tonci}.

The subalgebra of $\g$ corresponding to $Z(S)^0$ is precisely $\g_0$.

\begin{lemma}
$M_\Z$ acts cocompactly and properly discontinuously on a contractible
manifold $X_M$ and there is a sphere $\text{``}S_M\subset\partial M_\Z\text{''}$ with $\dim
S_M=\dim X_M-1$.
\end{lemma}

\begin{proof}
  $M_\Z$ is a lattice in $M_\R$ by \cite[Theorem
  9.4]{bhc:arithmetic}. It is a cocompact lattice since the $\Q$-rank
  of $M$ is 0. The manifold $X_M$ can be taken to be the product of
  the symmetric space of the semi-simple factor of $M$ and of the
  Euclidean factor corresponding to $T_\R/K_T$ (real points in the
  torus modulo maximal compact subgroup). See 
  Example \ref{uniform}.
\end{proof}

The Lie algebra of a minimal $\Q$-parabolic subgroup
$P$ is $\g_0\oplus \bigoplus_{\alpha\in\Phi^+}\g_\alpha$ and $P=Z(S)^0\cdot U$
where $U$ is the connected nilpotent subgroup with Lie algebra
$\bigoplus_{\alpha\in\Phi^+}\g_\alpha$. 

\begin{lemma}\label{parabolic}
  $P_\Z=P\cap GL_N(\Z)$ acts cocompactly and properly discontinuously
  on a contractible manifold $X_P$ and there is a sphere
  $\text{``}S_P\subset\partial P_\Z\text{''}$ with $\dim S_P=\dim
  X_P-1$ and $\dim X_P=\dim U_\R+\dim X_M$. Moreover, $\dim X_P=\dim
  G_\R/K-r$.
\end{lemma}

\begin{proof}
Note that $S$ intersects $G_\Z$ in a finite subgroup. After passing to
a finite cover, we have a split exact sequence
$$1\to U\to P/S\to Z(S)^0/S\to 1$$ and $Z(S)^0/S$ is a quotient of $M$
with finite kernel. The image of $P_\Z$ is a lattice in $P_\R/S_\R$. A
maximal compact subgroup $K'$ of $Z(S)_\R^0/S_\R$ lifts to a maximal
compact subgroup of $P_\R/S_\R$.  We set $X_P=(P_\R/S_\R)/K'$ and apply
Lemma \ref{ses} to $K'\rtimes U_\R\subset P_\R/S_\R$. To prove the last
statement, recall \cite{bs:corners} that $G_\Z$ acts cocompactly and
properly discontinuously on a contractible manifold with corners whose
interior can be identified with $G_\R/K$ and $P_\Z$ is the stabilizer
of a lowest dimensional stratum, which is a contractible manifold and
has codimension $r$ (and is really a copy of $X_P$).
\end{proof}

We now note that if the $\Q$-rank $r=0$ the theorem follows from Example
\ref{uniform}. From now on we will assume that $r\geq
1$. 

We now state some lemmas. For every $\alpha\in\Delta$ set $\hat\a=2\alpha$ if
$2\a\in \Phi$ and otherwise set $\hat\a=\a$. Let
$\hat\Delta=\{\hat\a|\a\in\Delta\}$.

\begin{lemma}\label{key}
  Let $\Phi$ be a (possibly unreduced,
  possibly not irreducible) root system, and let $\Delta$ be the set
  of simple roots (with respect to some ordering). There is
  an ordering of the set $\hat\Delta$ with the following
  properties. Let $\hat\a\in\hat\Delta$. Suppose that the elements of
  $\hat\Delta$ that precede $\hat\a$ in the order are labeled by one
  of the letters ``U'' or ``D'', and also label $\hat\a$ itself by
  ``D''. We refer to the elements of $\hat\Delta$ labeled ``D'' as
  {\it D-nodes}, and those labeled ``U'' as {\it U-nodes}. Then there
  exist $\sigma,\mu\in \Phi\cup\{0\}$ such that:
\begin{enumerate}
\item $\mu-\sigma=\hat\a$,
\item the difference $\sigma-\phi$ between $\sigma$ and any $\phi\in \Phi\cup\{0\}$ is
  not a positive multiple of a D-node,
\item the difference $\phi-\sigma$ between any $\phi\in \Phi\cup\{0\}$ is
  not a positive multiple of a U-node.
\end{enumerate}
\end{lemma}

Label ``U'' means that it is not possible to go ``up'' from $\sigma$ along the
simple root and reach a root or 0 and, similarly, ``D'' stands for
``down''. We postpone the proof of this lemma until the end of Section
\ref{sec:B}.

\begin{lemma}\label{fh}
\begin{enumerate}
\item If $v\in \g_\alpha$, $w\in\g_\beta$ then $[v,w]\in\g_{\a+\b}$ (the
  latter is defined to be 0 unless $\a+\b\in \Phi\cup\{0\}$).
\item Suppose that $\alpha,\beta,\alpha+\beta\in \Phi\cup\{0\}$ and $\a\neq 0$.
  Then there exist $v\in \g_\a$ and $w\in\g_\beta$ such that $[v,w]\neq 0$.
\end{enumerate}
\end{lemma}

\begin{proof} Both statements are well-known over $\C$ (see
  e.g. \cite[Proposition 2.5]{tonci} for (1) and \cite[Claim 21.19]{fh:reps}
  for (2)). (1) is proved the same way over $\Q$, and statement (2) over $\Q$
  follows by decomposing each root space into 1-dimensional subspaces which
  are root spaces with respect to a maximal torus that contains $S$ (but is
  split only over $\C$, not over $\Q$).\end{proof}

\begin{lemma}\label{exp}
$$\exp(X)Y\exp(-X)=Y+[X,Y]+\frac{1}{2!}[X,[X,Y]]+
\frac{1}{3!}[X,[X,[X,Y]]]+\cdots\qed$$\end{lemma}

\begin{lemma}{\cite[Corollary 10.14]{raghunathan:book}}\label{nocontraction}
  Let $M$ be an algebraic group defined over $\Q$ and let $\rho:M\to GL(V)$ be
  a homomorphism defined over $\Q$ into the general linear group of a vector
  space defined over $\Q$. If $\Gamma\subset M$ is an arithmetic lattice in
  $M$ and if $L$ is a lattice in $V_\Q$ then there is a finite index subgroup
  $\Gamma'$ of $\Gamma$ such that $\rho(\Gamma')$ preserves $L$. \qed
\end{lemma}

For each $\hat\a\in\hat\Delta$ choose a rational vector
$\zeta_{\hat\a}\in\g_{-\hat\a}$ such that
$$[\zeta_{\hat\a},\cdot]:\g_{\mu}\to\g_{\sigma}$$
is nonzero whenever
$\sigma,\mu\in\Phi\cup\{0\}$ are such that $\mu-\sigma=\ha$. That such a
vector exists follows from Lemma \ref{fh} and it can be taken to be a rational
vector by perturbing (the set of bad choices is contained in a finite union of
proper subspaces).

Let $\alpha_1,\cdots,\alpha_r$ be the simple roots in $\Delta$ ordered so that
$\hat\a_1,\hat\a_2,\cdots,\hat\a_r$ is the ordering of $\hat\Delta$ from Lemma
\ref{key}. Every root in $\Phi^+$ is an integral combination with nonnegative
coefficients of the simple roots.  For $i=1,\cdots,r$ denote by $\Phi_i^+$ the
set of positive roots obtained as nonnegative integral linear combinations of
simple roots $\alpha_1,\alpha_2,\cdots,\alpha_i$ and involving $\a_i$ with a
positive coefficient. Note that if $\alpha,\alpha'\in \Phi_i^+$ and if
$\alpha+\a'$ is a root, then $\alpha+\a'\in\Phi_i^+$. It follows that
$\n_i=\bigoplus_{\a\in\Phi_i^+}\g_\a$ is a nilpotent subalgebra of $\g$. The
integral subgroup $N_i=exp(\n_i)$ is nilpotent and the intersection $N_i\cap
G_\Z$ is a lattice in $N_i\cap G_\R$. Thus by Corollary \ref{nilpotent} we
have an expanding map $C_i\to N_i\cap G_\R$ from the cone $C_i$ on a sphere
$S_i$ with $$\dim C_i=\dim N_i\cap G_\R=\bigsum_{\a\in\Phi_i^+}{\dim \g_\a}$$
where $\dim\g_\a$ is the complex dimension (or equivalently the real dimension
of the real points of $\g_\a$). The $N_i$'s play the role of the ``column
groups'' in our examples. The reader will note that the maps $C_i\to N_i\cap
G_\R$ have been constructed as homeomorphisms and it will do no harm (and it
will simplify notation) to omit the name of the map and simply identify $C_i$
with $N_i$.

Similarly, for $i=1,2,\cdots,r$, consider the abelian subalgebra
$\g_{-\hat\a_i}$. By $P_i=exp(\g_{-\hat\a_i})$ denote the associated integral
subgroup and again by Corollary \ref{nilpotent} we know that $P_i\cap G_\Z$ is
a cocompact lattice in $P_i\cap G_\R$. The latter is a noncompact group and we
define a proper embedding $R_i:[0,\infty)\to P_i\cap G_\R$ (whose image we
also denote by $R_i$) by $$R_i(t)=exp(t\zeta_{\hat\alpha_i})$$ The
$P_i$'s play the role of the subdiagonal positions in our examples.

We now define
$$C=C_M\times (C_1\cup R_1)\times (C_2\cup
R_2)\times\cdots\times (C_r\cup R_r).$$ This is the open cone on a finite
complex $L$ obtained by taking joins of a sphere (for $C_M$) and spheres with a
point added. We take the join triangulation on $L$.  Define
also a map $$\Psi:C\to G_\R$$
as follows. Let $(d,z_1,z_2,\cdots,z_r)\in C$.
Then each $z_i$ is either equal to some $x_i\in C_i$ or to some $y_i\in R_i$.
Say $z_i=x_i\in C_i$ for $i=i_1,i_2,\cdots,i_k$ and $z_i=y_i\in R_i$ for
$i=j_1,j_2,\cdots,j_l$ where we assume $i_1<i_2<\cdots<i_k$,
$j_1<j_2<\cdots<j_l$ and $k+l=r$. Then set
$$\Psi(d,z_1,z_2,\cdots,z_r)=x_{i_1}x_{i_2}\cdots x_{i_k}dy_{j_1}y_{j_2}\cdots
y_{j_l}.$$ Note that the image of $\Psi$ is contained in a bounded neighborhood
of $G_\Z$ in $G_\R$ -- this is true for the components $R_i,C_i$ and
$C_M$ by construction and remains true after taking pointwise products
of such neighborhoods.  

\begin{remark}
Suppose $A$ and $B$ are $3\times 3$ elementary matrices with nonzero entries
in positions (1,2) and (2,3) respectively. Then $AB$ is a matrix with a
nonzero (1,3) position, while $BA$ is obtained from $A$ and $B$ by
``superposition'' (as in our examples). This explains why we have to
carefully arrange different components of the map $\Psi$.
\end{remark}

The proof of Theorem \ref{arithmetic} is now reduced to
the following two lemmas.
\vskip 0.5cm
\noindent
{\bf Lemma A.} {\it $L$ is an $m$-obstructor complex with $m=\dim
  X_G-2$.}\newline 
\vskip 0.01cm
\noindent
{\bf Lemma B.} {\it $\Psi$ is proper and expanding} 
\vskip 0.3cm

\begin{proof}[Proof of Lemma A]
  
  We have $$L=L_M*L_1*L_2*\cdots *L_r$$
  where $L_M$ is the sphere of dimension
  $\dim X_M-1$ and $L_i$ is the disjoint union of a sphere of dimension
  $\bigsum_{\a\in\Phi_i^+}{\dim \g_\a}-1$ and a point. Thus $L_M$ is a $(\dim
  X_M-2)$-obstructor complex and $L_i$ is a $\bigsum_{\a\in\Phi_i^+}{\dim
    \g_\a}-1$-obstructor complex.  The Join Lemma implies that $L$ is an
  $m$-obstructor complex for
\begin{equation*}
\begin{split}
  m=(\dim X_M-2)+\sum_{i=1}^r(\bigsum_{\a\in\Phi_i^+}{\dim
    \g_\a}-1)+2r=\\
  \dim X_M+\sum_{\a\in\Phi^+} 
{\dim \g_\a} +r-2=\dim X_P+r-2=\dim G_\R/K-2
\end{split}
\end{equation*}
where the last two equalities follow from Lemma \ref{parabolic} and
the observation that $\dim U_\R=\bigsum_{i=1}^r\bigsum_{\a\in\Phi_i^+}{\dim
    \g_\a}$.
\end{proof}

\section{Proof of Lemma B}\label{sec:B}

Choose two disjoint simplices of $L$ and a sequence of points in their
cones.  To simplify notation, we omit the subscripts corresponding to
the sequence counters. We then have $(d,z_1,\cdots,z_r)$ and
$(d',z_1',\cdots,z_r')$. Say the indices of $z_i$ corresponding to
nontrivial points in the rays are $j_1<j_2<\cdots<j_l$ and of $z_i'$
they are $k_1<k_2<\cdots<k_q$. The two sets are disjoint. For
concreteness, assume $j_l>k_q$. If possible, pass to a subsequence
such that the point $z_{j_l}=y_{j_l}\in R_{j_l}$ stays bounded. In
this case we replace the point $(d,z_1,\cdots,z_r)$ by the point in
which the $z_{j_l}$ coordinate is replaced by 1. This results in a
simpler pair of sequences and their divergence is equivalent to the
divergence of the original pair. We may thus assume that either
$y_{j_l}$ goes to infinity or that all $y_j$ and $y_j'$-coordinates
are 1.

First consider the case when $y_{j_l}\to \infty$. We will argue that the two
sequences diverge. It will not be important that the $x$-coordinates on the
two sides belong to disjoint simplices, so we will collect all $x$'s and
$x'$'s into an element denoted $u$, and we will let $s$ be the difference
between the two elements, so we write:
\begin{equation}
udy_{j_1}\cdots y_{j_l}=d'y_{k_1}'\cdots y_{k_q}'s
\end{equation}
where we assume, by way
of contradiction, that the sequence $s$ is bounded.

Consider the adjoint representation of $G$ on $Aut(\g)$. We will obtain a
contradiction by comparing the automorphisms of $\g$ the two sides
induce. 

\def\ha{\hat\alpha} We label the roots
$\hat\alpha_{j_1},\ha_{j_2},\cdots,\ha_{j_l}$ by ``D'' and
$\hat\alpha_{k_1},\ha_{k_2},\cdots,\ha_{k_q}$ by ``U''. Lemma \ref{key}
(with $\ha=\ha_{j_l}$) provides us with $\sigma,\mu\in \Phi\cup\{0\}$. We
will focus on the $\g_\mu\to\g_\sigma$ component of the transformation
(2). This component will play the role of a ``matrix position'' in our
examples.

First, we look at $udy_{j_1}\cdots y_{j_l}$. Lemma \ref{fh}(1) and
Lemma \ref{exp} imply that $Ad(y_{j_l})$ maps $\g_\mu$ into
$$\g_\mu\oplus
\g_{\mu-\ha_{j_l}}\oplus\g_{\mu-2\ha_{j_l}}\oplus\cdots$$ which, by
Lemma \ref{key}(1) and (2) is just $\g_\mu\oplus\g_\sigma$, and
further, $Ad(y_{j_l})=I\oplus [t\zeta_{\ha_{j_l}},\cdot]$ on
$\g_\mu$. By Lemma \ref{key}(2) we have $$Ad(y_{j_1}y_{j_2}\cdots
y_{j_{l-1}})|\g_\sigma =I$$ and 
$$Ad(y_{j_1}y_{j_2}\cdots
y_{j_{l-1}})|\g_\mu$$ has its image contained in 
\begin{equation}
\oplus\{\g_\phi|\mu-\phi\text{ is a nonnegative linear combination of
}\ha_{j_1},\cdots,\ha_{j_{l-1}}\}.
\end{equation}
Next, $Ad(d)$ preserves each $\g_\phi$ and $Ad(u)$ maps $\g_\phi$ into
$$\bigoplus_{\nu-\phi\in\Phi^+\cup\{0\}}\g_\nu.$$ We now claim that $Ad(u)$ has
$\g_\sigma$-component equal to 0 when restricted to (3). Otherwise, we
get equations:
\begin{equation*}
\begin{split}
\mu-\phi=&a_1\ha_{j_1}+\cdots+a_{l-1}\ha_{j_{l-1}}\\
\sigma-\phi=&b_1\ha_1+\cdots+b_r\ha_r
\end{split}
\end{equation*}
with all coefficients $\geq 0$
which, together with $\mu-\sigma=\ha_{j_l}$, imply (by subtracting the
second and third from the first) that
$$b_1\ha_1+\cdots+b_r\ha_r=a_1\ha_{j_1}+\cdots+a_{l-1}\ha_{j_{l-1}}-\ha_{j_l}$$
thus violating the linear independence of the simple roots.

Summarizing, we have

\begin{lemma}\label{lhs}
The $\g_\mu\to\g_\sigma$-component of $Ad(udy_{j_1}\cdots y_{j_l})$ is
$$D_\sigma\circ [t\zeta_{\ha_{j_l}},\cdot]$$ where $D_\sigma$ is the
restriction of $Ad(d)$ to $\g_\sigma$.
\qed\end{lemma}

We now perform a similar analysis for the right hand side of (2).

We have that $Ad(d')$ preserves each root space. By Lemma \ref{key}
(3), the only root space $\g_\phi$ such that $Ad(y_{k_1}'\cdots
y_{k_q}')$ has a nontrivial $\g_\sigma$-component when restricted to
it is $\g_\sigma$ and the component is $I$. We therefore have

\begin{lemma}\label{rhs}
The $\g_\mu\to\g_\sigma$-component of $Ad(d'y_{k_1}'\cdots y_{k_q}'s)$ is
$$D'_\sigma S_{\mu\sigma}$$ where $D'_\sigma$ is the restriction of
$Ad(d')$ to $\g_\sigma$ and $S_{\mu\sigma}$ is the
$\g_\mu\to\g_\sigma$-component of $Ad(s)$.\qed
\end{lemma}

We now conclude that $$D_\sigma\circ [t\zeta_{\ha_{j_l}},\cdot]=D'_\sigma
S_{\mu\sigma}$$
i.e. $$D\circ [t\zeta_{\ha_{j_l}},\cdot]$$
is a bounded sequence of maps 
as $t\to\infty$ and $D$ represents a sequence in $M_\R$. Let $L$ be a
(rational) lattice in $\g_\sigma$ that contains a nonzero vector $h$ in the
image of the map $[\zeta_{\ha_{j_l}},\cdot]$. Lemma \ref{nocontraction}
yields a subgroup $\Gamma$ of finite index in $M_\Z$ such that $Ad(\Gamma)$
preserves $L$. Since $\Gamma$ is cocompact in $M$, it follows that $Ad(M)$
{\it nearly preserves} $L$ in the sense that there is a
neighborhood $\Omega$ of $0\in\g_\sigma$ such that $Ad(M)$ sends no nonzero
element of $L$ (in particular, $h$) into $\Omega$ (see Mahler's criterion,
\cite[Corollary 10.9]{raghunathan:book}). It follows that $Ad(M)$ sends $th$
outside $t\Omega$, so the composition $D\circ [t\zeta_{\ha_{j_l}},\cdot]$ is
not a uniformly bounded map as $t\to\infty$. This contradiction proves Lemma
B under our original assumption that $y_{j_l}\to\infty$.

Now suppose that all $y_j$ and $y_j'$ coordinates are 1. Equation (2) becomes
\begin{equation}
ud=d's
\end{equation}
This is an equation in the minimal parabolic subgroup $P=Z(S)^0\cdot
U$. After passing to a finite cover, we may assume that
$P=Z(S)^0\rtimes U$. Consider the canonical homomorphism $\tau:P\to
Z(S)^0$. It maps $ud$ to $d$ and $d's$ to $d's'$ where $s'$ is a
bounded sequence in $M_\R$. Thus $d$ and $d'$ stay a bounded distance
apart. Since by construction they belong to divergent cones in $M_\R$,
they both have to stay bounded. It follows that $u$ must stay bounded,
i.e. that the sequence $x$ and the sequence $x'$ stay within bounded
distance. But by Corollary \ref{nilpotent} this means that both $x$ and
$x'$ stay bounded, contradicting the assumption that the original
sequences were chosen to be unbounded. This contradiction proves Lemma B.

\begin{proof}[Proof of Lemma \ref{key}]
Every root can be written as a linear
combination of simple roots. The key to this proof is the fact that
the coefficients in such linear combinations are explicitly known (and can be
found e.g. in \cite[Appendix C]{tonci}). 

We first note that if $\Phi$ is not irreducible, then the statement follows
immediately from the corresponding statements for the irreducible components.
The reason for this is that $\Phi$ is then the disjoint union of its
irreducible components, and any ordering of $\hat\Delta$ that restricts
correctly to each component will work. Moreover, if there are elements of
$\hat\Delta$ that are not labeled at all (which is the case when $\ha$ is not
the highest node) then we may restrict our consideration to the root system
$\Phi'$ generated by the labeled nodes (and in fact to the irreducible
component of $\Phi'$ that contains $\ha$).

Here is one situation when we can take $\sigma$ to be the negative of the
largest positive root: $\sigma+\ha\in \Phi\cup\{0\}$ and if we write
$\sigma=-m\ha-\cdots$ as a linear combination of simple roots, then all other
roots in $\Phi$ have their $\ha$-coefficient $>-m$. These situations, together
with the coefficients, are pictured below, with the node $\ha$ circled.
\begin{figure}[htbp]
\begin{center}
\input{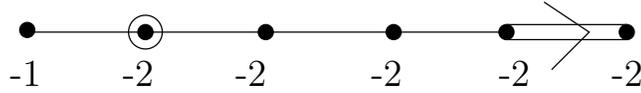}
\caption{\bf $B_n$, $n\geq 3$}
\end{center}
\end{figure}
\begin{figure}[htbp]
\begin{center}
\input{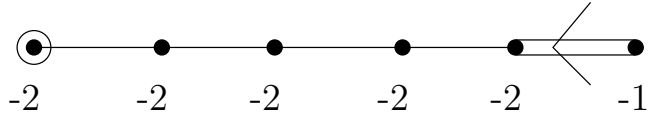}
\caption{\bf $C_n$, $n\geq 2$}
\end{center}
\end{figure}
\begin{figure}[htbp]
\begin{center}
\input{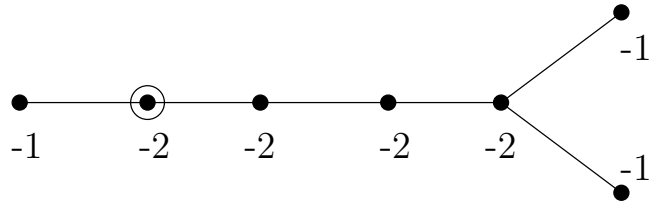}
\caption{\bf $D_n$, $n\geq 4$}
\end{center}
\end{figure}
\begin{figure}[htbp]
\begin{center}
\input{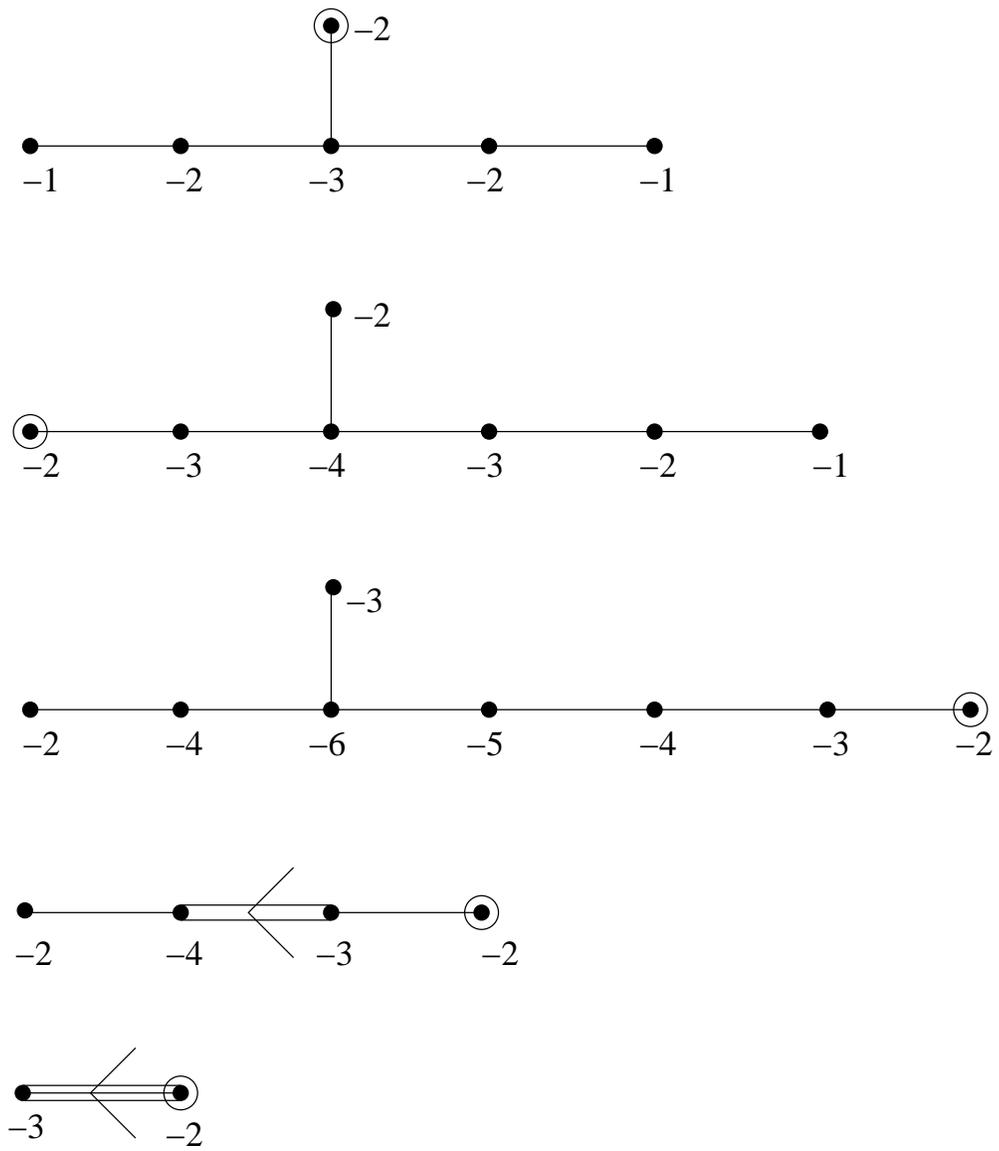}
\caption{\bf Exceptional reduced root systems $E_6,E_7,E_8,F_4,G_2$}
\end{center}
\end{figure}
\begin{figure}[htbp]
\begin{center}
\input{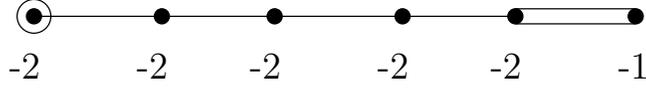}
\caption{\bf Unreduced root system $BC_n$, $n\geq
2$. The positive roots are $e_i\pm e_j$ for $i<j$, $e_i$, $2e_i$,
$i=1,2,\cdots,n$. The nodes correspond to the roots (from left to right)
$e_1-e_2,e_2-e_3,\cdots,e_{n-1}-e_n,e_n,2e_n$ with the last two
corresponding to the same (rightmost) node in the diagram. The largest
root is $2e_1$ and its coefficients with respect to the nodes are
$2,2,\cdots,2,1$ if we take $\ha_n=2e_n$ as the representative of the
last node.}
\end{center}
\end{figure}

In the case of $BC_1$ (i.e. $\Phi=\{-2\a,-\a,\a,2\a\}$) we can take
$\sigma=-2\a=-\ha$ and $\mu=0$. This leaves us with type $A_n$. The ordering
is the usual linear ordering of the nodes. Say $\a$ is node $k$ in this
ordering and let a labeling by ``U'''s and ``D'''s of nodes $\leq k$ be given
as in the lemma. Let $1\leq l\leq k$ be such that the label of node $l$ is
``D'' but the label of node $l-1$ is not ``D'' (or $l=1$). Define $\sigma$ as
the negative of the sum of the nodes $l,l+1,\cdots,k$.

In general, to define the ordering on $\hat\Delta$, we follow this procedure:
Work separately on components. If a component is of type $A_n$ order the nodes
linearly. Otherwise, define the highest root in the ordering to be the circled
node in the corresponding figure above. Then pass to the subdiagram consisting
of unlabeled nodes and repeat the procedure. 
\end{proof}

We finish the paper by looking at few more examples.

\begin{remark}
  The ordering of $\Delta$ was used in two places: to define $\Phi_i^+$ (and
  the associated nilpotent groups $N_i$) and in determining the order in which
  the $y_i$'s come in the definition of the map $\Psi$. The reader will
  observe that we can use two different orders on $\Delta$: an arbitrary order
  to define $\Phi_i^+$ and the one from Lemma \ref{key} to order the $y_i$'s.
  It is convenient to use the standard order for the first purpose since then
  the $N_i$'s are standard nilpotent matrix groups. This is the practice we
  follow in the examples.
\end{remark}

\section{$\G=Sp_{2n}(\Z)$}

Consider $V=\C^{2n}$ with the standard $\C$-basis
$e_1,e_2,\cdots,e_n,e_{n+1},\cdots,e_{2n}$. Let $J$ be the
anti-symmetric bilinear pairing defined by $J(e_i,e_{n+i})=1$ and
$J(e_i,e_j)=0$ if $|j-i|\neq n$. The group $Sp_{2n}(\C)$ is the subgroup
of $GL_{2n}(\C)$ consisting of matrices that preserve $J$. It is
convenient to represent the matrices in $Sp_{2n}(\C)$ and in its Lie
algebra in $2\times 2$ block form corresponding to the partition of
the basis for $V$ into the first $n$ and last $n$ vectors. Thus $J$ is
represented by $J=\begin{pmatrix} 0&-I\\ I&0\end{pmatrix}$. The Lie algebra
consists of block matrices $X=\begin{pmatrix} A&B\\ C&D\end{pmatrix}$ such that
$XJ+JX^t=0$, i.e. $B$ and $C$ are symmetric matrices, and $D$ is the
negative transpose of $A$. A maximal split torus $S$ can be taken to
consist of diagonal matrices
\begin{equation}
diag[s_1,s_2,\cdots,s_n,s_1^{-1},s_2^{-1},\cdots,s_n^{-1}]
\end{equation}
Denote
by $y_i:S\to\C$ the character that takes the above matrix to
$s_i$. The positive roots are $$y_i-y_j, (i<j); y_i+y_j, (i\neq j);
2y_i$$ and the simple roots are
$y_1-y_2,y_2-y_3,\cdots,y_{n-1}-y_n,2y_n$, so the root system is of
type $C_n$. The Lie algebra $\frak n_i$ for $i<n$ consists of matrices
that in the column $i+1$ of the $(1,1)$-block have arbitrary entries
above the diagonal and 0 in all other positions, the $(2,2)$-block is
the negative transpose of the $(1,1)$-block, and the $(1,2)$ and
$(2,1)$-blocks are 0. For $i=n$ the Lie algebra $\frak n_n$ consists
of matrices $\begin{pmatrix} 0&B\\ 0&0\end{pmatrix}$ with $B$ symmetric. In
this case $\Phi_n^+$ consists of the roots in classes 2 and 3
above. The obstructor complex in this case is
$$L=S^0_+*S^1_+*\cdots *S^{n-2}_+*S^{\frac{(n+2)(n-1)}{2}}_+$$ The Lie
algebras $\frak n_i$ are abelian and exponentiation (to obtain $N_i$)
amounts to adding the identity matrix.

\section{$\G=Sp_{2n}(\O)$}

If $r$ and $s$ denote the numbers of real and complex places, then
$\O\cong \Z^{r+2s}$ as an abelian group so the dimensions of $N_i$'s
above should be multiplied by $r+2s$. The centralizer $Z(S)$ consists
of the diagonal matrices (5) and $M_\Z$ is commensurable with the
group of such matrices with entries in $\O^*$: it is an abelian group
of rank $(n-1)(r+2s)$. The obstructor complex is thus
$$
S^{(n-1)(r+2s)-1}*S^{r+2s-1}_+*S^{2(r+2s)-1}_+*\cdots*S^{(n-1)(r+2s)-1}_+*S^{\frac{n(n+1)(r+2s)}{2}-1}_+$$ 

\section{$\G=SO(Q)$ for a nondegenerate form $Q$}

Any nondegenerate quadratic form $Q$ defined over $\Q$ on a vector
space $V$ can be represented as a direct sum of a certain number, say
$q$, of hyperbolic planes and of an anisotropic quadratic form $Q_0$
(i.e. $Q_0$ does not take value 0 on nonzero rational vectors). We
follow notation from \cite[V.23.4]{borel:lag}, where the reader can
find more details about $SO(Q)$. Choose a rational basis
$e_1,e_2,\cdots,e_n$ of $V$ so that $e_i,e_{n-q+i}$ span a hyperbolic
direct summand for $i=1,2,\cdots,q$ (with
$<e_i,e_i>=<e_{n-q+i},e_{n-q+i}>=0, <e_i,e_{n-q+i}>=1$ in the
associated symmetric pairing) and $e_{q+1},\cdots,e_{n-q}$ span
the $Q_0$-summand. The maximal split $\Q$-torus $S$ in $SO(Q)$ is the
group of diagonal matrices
$$diag[s_1,s_2,\cdots,s_q,1,1,\cdots,1,s_1^{-1},s_2^{-1},\cdots,s_q^{-1}]$$
Denote by $y_i:S\to\C$ the character that sends the above matrix to
$s_i$. The centralizer $Z(S)$ is $S\times SO(Q_0)$ and the positive
roots are $$y_i-y_j, (i<j); y_i+y_j, (i\neq j); y_i $$ with the first
two kinds having multiplicity 1 and the third kind of multiplicity
$n-2q$ (so the third kind is not present when $n=2q$). The simple
roots are $y_1-y_2,y_2-y_3,\cdots,y_{q-1}-y_q,y_q$ which is type $B_q$
(if $n=2q$ the type is $D_q$). Every element of $SO(Q)$ and of its Lie
algebra is conveniently represented in a $(3,3)$-block form
corresponding to first $q$, middle $n-2q$, and last $q$ basis
vectors. A block matrix $(A_{ij})$ is in the Lie algebra of $SO(Q)$ if
and only if the following 6 conditions are satisfied:
\begin{equation*}
\begin{split} 
A_{13}+A_{13}^t=A_{12}F_0+A_{23}^t=A_{11}+A_{33}^t=0\\
A_{22}F_0+F_0A_{22}^t=A_{32}F_0+A_{21}^t=A_{31}+A_{31}^t=0
\end{split}
\end{equation*} 
where $F_0$ is the matrix of $Q_0$.

We now describe the Lie algebras ${\frak n}_i$ of the groups $N_i$,
$i=1,2,\cdots,q-1$ (we are
using the standard ordering of the simple roots): In the $(1,1)$-block
consider the $i$ positions above the diagonal in column $i+1$. Any
(complex) entries are allowed. Then change the sign of these entries
and transpose, and write this row vector in the $3\times 3$ block in
row $i+1$ to the left of the diagonal. To get $N_i$ exponentiate -- in
this case this amounts to adding the identity matrix. Intersecting
with real points, this defines the space $C_i$ from the proof of
Lemma B  -- it is the cone on the sphere of dimension $i-1$.

The group $N_q$ is not abelian -- it is 2-step nilpotent. Its Lie
algebra is spanned by the positive roots in classes 2 and 3 above and
it consists of block matrices that have vanishing all 3 diagonal
blocks and all 3 blocks below the diagonal. The group of real points
of $N_q$ has dimension $q(n-2q)+(1+2+\cdots+(q-1))$.

The ray corresponding to the root $y_i-y_{i+1}$ is obtained by
exponentiating matrices with entry $t\geq 0$ in position $(i+1,i)$ of
$A_{11}$ and entry $-t$ in position $(i,i+1)$ of $A_{33}$ -- in this case
exponentiation amount to adding the identity matrix. To define a ray
corresponding to the root $y_q$ choose a nonzero rational row vector
$v$ of length $n-2q$ and place it in the last row of $A_{12}$ and then
place the column vector $-F_0 v^t$ in the last column of $A_{23}$. The
ray is obtained by exponentiating positive multiples of this matrix --
this amounts to adding the identity matrix and the entry
$-\frac{1}{2}vF_0v^t$ in position $(q,q)$ of $A_{13}$.

If the form $Q_0$ is definite, then $M_\R=SO(Q_0)_\R$ is a compact
group and $C_M$ will be a point. If $Q_0$ is not definite, then $C_M$
can be identified with the symmetric space of $SO(Q_0)_\R$.

\bibliography{$HOME/lib/ref} 

\providecommand{\bysame}{\leavevmode\hbox to3em{\hrulefill}\thinspace}
\begin{thebibliography}{BHC62}

\bibitem[Ash77]{ash:spine}
Avner Ash, \emph{Deformation retracts with lowest possible dimension of
  arithmetic quotients of self-adjoint homogeneous cones}, Math. Ann.
  \textbf{225} (1977), no.~1, 69--76.

\bibitem[BH99]{bh:book}
Martin~R. Bridson and Andr{\'e} Haefliger, \emph{Metric spaces of non-positive
  curvature}, Springer-Verlag, Berlin, 1999.

\bibitem[BHC62]{bhc:arithmetic}
Armand Borel and Harish-Chandra, \emph{Arithmetic subgroups of algebraic
  groups}, Ann. of Math. (2) \textbf{75} (1962), 485--535.

\bibitem[BKK]{bkk:vk}
Mladen Bestvina, Michael Kapovich, and Bruce Kleiner, \emph{Van {K}ampen's
  embedding obstruction for discrete groups}, preprint, 2000.

\bibitem[Bor69]{borel:ga}
Armand Borel, \emph{Introduction aux groupes arithm\'etiques}, Hermann, Paris,
  1969, Publications de l'Institut de Math\'ematique de l'Universit\'e de
  Strasbourg, XV. Actualit\'es Scientifiques et Industrielles, No. 1341.

\bibitem[Bor91]{borel:lag}
Armand Borel, \emph{Linear algebraic groups}, second ed., Springer-Verlag, New
  York, 1991.

\bibitem[BS73]{bs:corners}
A.~Borel and J.-P. Serre, \emph{Corners and arithmetic groups}, Comment. Math.
  Helv. \textbf{48} (1973), 436--491, Avec un appendice: Arrondissement des
  vari\'et\'es \`a coins, par A. Douady et L. H\'erault.

\bibitem[BW80]{bw:book}
Armand Borel and Nolan~R. Wallach, \emph{Continuous cohomology, discrete
  subgroups, and representations of reductive groups}, Princeton University
  Press, Princeton, N.J., 1980.

\bibitem[FH91]{fh:reps}
William Fulton and Joe Harris, \emph{Representation theory}, Springer-Verlag,
  New York, 1991, A first course, Readings in Mathematics.

\bibitem[Got48]{goto}
Morikuni Got\^o, \emph{Faithful representations of {L}ie groups. {I}}, Math.
  Japonicae \textbf{1} (1948), 107--119.

\bibitem[GS92]{gs:harmonic}
Mikhail Gromov and Richard Schoen, \emph{Harmonic maps into singular spaces and
  $p$-adic superrigidity for lattices in groups of rank one}, Inst. Hautes
  \'Etudes Sci. Publ. Math. (1992), no.~76, 165--246.

\bibitem[Hel78]{helgason:book}
Sigurdur Helgason, \emph{Differential geometry, {L}ie groups, and symmetric
  spaces}, Academic Press Inc. [Harcourt Brace Jovanovich Publishers], New
  York, 1978.

\bibitem[Hoc65]{hochschild}
G.~Hochschild, \emph{The structure of {L}ie groups}, Holden-Day Inc., San
  Francisco, 1965.

\bibitem[Iwa49]{iwasawa:k}
Kenkichi Iwasawa, \emph{On some types of topological groups}, Ann. of Math. (2)
  \textbf{50} (1949), 507--558.

\bibitem[Kna96]{tonci}
Anthony~W. Knapp, \emph{Lie groups beyond an introduction}, Birkh\"auser Boston
  Inc., Boston, MA, 1996.

\bibitem[Mal45]{malcev:k}
A.~Malcev, \emph{On the theory of the {L}ie groups in the large}, Rec. Math.
  [Mat. Sbornik] N. S. \textbf{16(58)} (1945), 163--190.

\bibitem[Mil76]{millson:betti}
John~J. Millson, \emph{On the first {B}etti number of a constant negatively
  curved manifold}, Ann. of Math. (2) \textbf{104} (1976), no.~2, 235--247.

\bibitem[Mos55]{mostow:sag}
G.~D. Mostow, \emph{Self-adjoint groups}, Ann. of Math. (2) \textbf{62} (1955),
  44--55.

\bibitem[Rag72]{raghunathan:book}
M.~S. Raghunathan, \emph{Discrete subgroups of {L}ie groups}, Springer-Verlag,
  New York, 1972, Ergebnisse der Mathematik und ihrer Grenzgebiete, Band 68.

\bibitem[R{\u{S}}97]{repovs:hilbert}
Du{\u{s}}an Repov{\u{s}} and Evgenij {\u{S}}{\u{c}}epin, \emph{A proof of the
  {H}ilbert-{S}mith conjecture for actions by {L}ipschitz maps}, Math. Ann.
  \textbf{308} (1997), no.~2, 361--364.

\bibitem[Zim84]{zimmer:book}
Robert~J. Zimmer, \emph{Ergodic theory and semisimple groups}, Birkh\"auser
  Verlag, Basel, 1984.

\end{thebibliography}
\end{document}